\documentclass[12pt]{article}
\usepackage[T2A,OT1]{fontenc}
\usepackage[ot2enc]{inputenc}
\usepackage[russian,english]{babel}
\usepackage{amsthm}
\def\C{\mathbf{C}}
\def\bC{\mathbf{\overline{C}}}
\newtheorem{theorem}{Theorem}
\newtheorem*{prop}{Proposition}
\newtheorem{lemma}{Lemma}
\theoremstyle{remark}
\newtheorem*{ack}{Acknowledgment}
\begin{document}
\author{Walter Bergweiler\thanks{Supported by the Deutsche
Forschungsgemeinschaft, Be 1508/7-1,
 and the ESF Networking Programme HCAA.} \ and
Alexandre Eremenko\thanks{Supported by NSF grant
DMS-1067886.}}
\title{A property of the derivative of an entire function}
\maketitle
\begin{abstract} We prove that the derivative of a 
non-linear entire function
is unbounded on the preimage of an unbounded set.

\end{abstract}

\section{Introduction and results}

The main result of this paper is the following theorem
conjectured by Allen Weitsman (private communication):
\begin{theorem} Let $f$ be a
non-linear entire function and $M$ an
unbounded set in $\C$. Then $f'(f^{-1}(M))$ is unbounded.
\end{theorem}

We note that there exist entire functions $f$ such that
$f'(f^{-1}(M))$ is bounded for every bounded set $M$,
for example, $f(z)=e^z$ or $f(z)=\cos z$.

Theorem 1 is a consequence of the following stronger result:
\begin{theorem}
Let $f$ be a transcendental entire function and $\varepsilon>0$.
Then there exists $R>0$ such that for every
$w\in\C$ satisfying $|w|>R$ there exists $z\in\C$ with 
$f(z)=w$ and 
$|f'(z)|\geq|w|^{1-\varepsilon}$. 
\end{theorem}

The example $f(z)=\sqrt{z}\sin\sqrt{z}$ shows that
that the exponent $1-\varepsilon$ in the last inequality cannot
be replaced by $1$. The function $f(z)=\cos \sqrt{z}$ has the property
that for every $w\in\C$ we have $f'(z)\to 0$ as $z\to\infty$, 
$z\in f^{-1}(w)$.

We note that the Wiman--Valiron theory 
\cite{V2,Hayman1974,BRS}
says that there exists a set $F\subset [1,\infty)$ of
finite logarithmic measure such that if
$$|z_r|=r\notin F
\quad\mbox{and}\quad
|f(z_r)|=\max_{|z|=r}|f(z)|,$$
then
$$
f(z)\sim \left(\frac{z}{z_r}\right)^{\nu(r,f)}f(z_r)
\quad\mbox{and}\quad
f'(z)\sim \frac{\nu(r,f)}{r}f(z)$$
for $|z-z_r|\leq r\nu(r,f)^{-1/2-\delta}$
as $r\to\infty$.
Here $\nu(r,f)$ denotes the central index and
$\delta>0$.
This implies that the conclusion of Theorem~2 holds for 
all $w$ satisfying $|w|=M(r,f)$ for some sufficiently
large $r\notin F$.
However, in general the exceptional set in 
the Wiman--Valiron theory is non-empty (see, e.g.,~\cite{B90})
and thus it seems that our results cannot be proved
using Wiman--Valiron theory.

\begin{ack}
We thank Allen Weitsman for helpful discussions.
\end{ack}
\section{Preliminary results} 

One important tool in the proof is the following result known as the Zalcman
Lemma \cite{Z}.
Let $$g^\#=\frac{|g'|}{1+|g|^2}$$
denote the spherical derivative of a meromorphic function~$g$.
\begin{lemma} Let $F$ be a non-normal family of meromorphic functions in
a region $D$. Then there exist a sequence $(f_n)$ in $F$, a sequence
$(z_n)$ in $D$, a sequence $(\rho_n)$ of positive real numbers and a
non-constant function $g$ meromorphic in $\C$ such that $\rho_n\to 0$
and $f_n(z_n+\rho_nz)\to g(z)$ locally uniformly in~$\C$.
Moreover, $g^\#(z)\leq g^\#(0)=1$ for $z\in\C$.
\end{lemma}

We say that $a\in\bC$ is a {\em totally ramified} value of a
meromorphic function $f$ if all $a$-points of $f$ are multiple.
A classical result of Nevanlinna says that a non-constant function meromorphic
in the plane can have at most $4$ totally ramified values, and that a
non-constant entire function can have at most $2$
finite totally ramified values.
Together with  Zalcman's Lemma this 
yields the following 
result~\cite{ChenHua,Hink,Lappan}; cf.~\cite[p.~219]{Zalcman98}.
\begin{lemma} Let $F$ be a family of functions meromorphic in a domain $D$
and $M$ a subset of $\bC$ with at least $5$ elements. Suppose that
there exists $K\geq 0$ such that for all $f\in F$ and $z\in D$ the
condition $f(z)\in M$ implies $|f'(z)|\leq K$. Then $F$ is a normal family.

If all functions in $F$ are holomorphic, then the conclusion holds 
if $M$ has at least $3$ elements.
\end{lemma}

Applying Lemma 2 to the family $\{ f(z+c):c\in\C\}$ where $f$ is an entire
function, we obtain the following result.

\begin{lemma} Let $f$ be an entire function and $M$ a subset of $\C$
with at least $3$ elements. If $f'$ is bounded on $f^{-1}(M)$, then
$f^\#$ is bounded in $\C$.
\end{lemma}

It follows from Lemma 3 that the conclusion of Theorems 1 and 2 holds for all
entire functions for which $f^\#$ is unbounded.

We thus consider entire functions with bounded spherical derivative.
The following result is due to Clunie and Hayman \cite{CH}. Let 
$$M(r,f)=\max_{|z|\leq r}|f(z)|
\quad\mbox{and}\quad 
\rho(f)=\limsup_{r\to\infty}\frac{\log\log M(r,f)}{\log r}$$
denote the maximum modulus and
the order of $f$.

\begin{lemma} Let $f$ be an entire function for which $f^\#$ is bounded.
Then $\log M(r,f)=O(r)$ as $r\to\infty$. In particular, $\rho(f)\leq 1$.
\end{lemma}
We will include a proof of Lemma 4 after Lemma~6.

The following result is due to Valiron \cite[III.10]{V2} and
H. Selberg \cite[Satz~II]{S}.

\begin{lemma} Let $f$ be a non-constant
 entire function of order at most $1$ for which
$1$ and $-1$ are totally ramified. Then $f(z)=\cos(az+b)$, where $a,b\in\C$,
$a\neq 0$.
\end{lemma}
We sketch the proof of Lemma 5.
Put
$h(z)=f'(z)^2/(f(z)^2-1)$. 
Then $h$ is entire and the lemma
on the logarithmic derivative~\cite[p.94, (1.17)]{Goldberg08},
together with the hypothesis
that $\rho(f)\leq 1$, yields that $m(r,h)=o(\log r)$ 
and hence that $h$ is constant.  This implies that 
$f$ has the form given.
Another proof is given in~\cite{GT}

The next lemma can be extracted from the work of Pommerenke
\cite[Sect. 5]{P}, see \cite[Theorem~5.2]{E}. 

\begin{lemma} Let $f$ be an entire function and $C>0$. If $|f'(z)|\leq C$
whenever $|f(z)|=1$, then $|f'(z)|\leq C|f(z)|$ whenever $|f(z)|\geq 1$.
\end{lemma}

Lemma 6 implies the theorem of Clunie and Hayman 
mentioned above (Lemma~4).
For the convenience of the reader 
we include a proof of a slightly more general statement, which is also more
elementary than the proofs of Clunie, Hayman and Pommerenke; see 
also~\cite[Lemma~1]{BarrettE}.

Let $G=\{ z:|f(z)|>1\}$ and $u=\log|f|$. Then $|f'/f|=|\nabla u|$ and
our statement which implies Lemmas 4 and 6 is the following.

\begin{prop}
Let $G$ be a region in the plane, $u$ a harmonic function in
$\overline{G}$, positive in
$G$, and such that for $z\in\partial G$ we have $u(z)=0$ and
$|\nabla u(z)|\leq 1$. 
Then $|\nabla u(z)|\leq 1$ for $z\in G$, and $u(z)\leq|z|+O(1)$
as $z\to\infty$.
\end{prop}

{\em Proof.} It is enough to consider the case of
unbounded $G$ with non-empty boundary.
For $a\in G$, consider the largest disc $B$ centered at $a$ and contained
in $G$. The radius 
$d=d(a)$ of this disc is the distance from $a$ to $\partial G$.
There is a point $z_1\in\partial B$ such that
$u(z_1)=0$. Put $z(r)=a+r(z_1-a),$ where $r\in(0,1)$.
Harnack's inequality gives
$$\frac{u(a)}{d(1+r)}\leq
\frac{u(z(r))}{d(1-r)}=\frac{u(z(r))-u(z_1)}{d(1-r)}.$$
Passing to the limit as $r\to 1$ we obtain 
$$u(a)\leq 2d(a)|\nabla u(z_1)|\leq 2d(a).$$
This holds for all $a\in G$. Now we take the gradient of both
sides of the Poisson formula and,
noting that 
$u(a+d(a)e^{it})\leq 2 d(a+d(a)e^{it})\leq 4 d(a)$,
 obtain the estimate
$$|\nabla u(a)|\leq\frac{1}{\pi d(a)}\int_{-\pi}^\pi
|u(a+d(a)e^{it})|dt\leq 8.$$
So $\nabla u$ is bounded in $G$. As  the complex conjugate
of $\nabla u$ is
holomorphic in $G$ and $|\nabla u(z)|\leq 1$ at all
boundary points $z$
of $G$, except infinity, the Phragm\'en--Lindel\"of theorem
\cite[III, 335]{PS} gives
that $|\nabla u(z)|\leq 1$ for $z\in G$. This completes the proof 
of the Proposition.

\medskip

We recall that for a non-constant entire function $f$ the maximum modulus
$M(r)=M(r,f)$ is a continuous strictly increasing function of $r$.
Denote by $\varphi$ the inverse function of~$M$.
Clearly, for $|w|>|f(0)|$ the equation
$f(z)=w$ has no solutions in the open disc of radius $\varphi(|w|)$ around $0$.
The following result of Valiron (\cite{V0,V}, see also \cite{EF}) says that for
functions of finite order this equation has solutions in a somewhat
larger disc.

\begin{lemma} Let $f$ be a transcendental entire function of finite order
and $\eta>0$. Then there exists $R>|f(0)|$ such that for all
$w\in\C,|w|\geq R$,
the equation $f(z)=w$ has a solution $z$
satisfying 
$|z|<\varphi(|w|)^{1+\eta}.$  
\end{lemma}

We note that Hayman (\cite{H}, see also \cite[Theorem~3]{B}) has constructed
examples which show that the assumption about finite order is essential
in this lemma.

\section{Proof of Theorem 2 }

Suppose that the conclusion is false. Then there exists $\varepsilon>0$,
a transcendental entire function $f$ and a sequence $(w_n)$ tending
to $\infty$
such that $|f'(z)|\leq|w_n|^{1-\varepsilon}$ whenever $f(z)=w_n$.
By Lemma 3, the spherical derivative of $f$ is bounded, and we may assume without loss of generality that
\begin{equation}
\label{3a}
f^\#(z)\leq 1\quad\mbox{for}\ z\in\C.
\end{equation}
We may also assume that $f(0)=0$. It follows from (\ref{3a})
that $|f'(z)|\leq 2$ if $|f(z)|=1$, and thus Lemma 6
yields
\begin{equation}\label{3b}
\left|\frac{f'(z)}{f(z)}\right|\leq 2\quad\mbox{if}\ |f(z)|\geq 1.
\end{equation}
It also follows from (\ref{3a}), together with Lemma 4, that $\rho(f)\leq 1$.
We may thus apply Lemma 7 and find that if $\eta>0$ and if $n$ is sufficiently
large, then there exists $\xi_n$ satisfying
$$|\xi_n|\leq \varphi(|w_n|)^{1+\eta}\quad\mbox{and}\quad f(\xi_n)=w_n.$$
We put
$$\tau_n=\varphi(|w_n|)^{1+2\eta}$$
and define
$$\Phi_n(z)=\frac{w_n-2f(\tau_nz)}{w_n}=1-2\frac{f(\tau_nz)}{w_n}.$$
Then $\Phi_n(0)=1,\;\Phi_n(\xi_n/\tau_n)=-1,$ and $\xi_n/\tau_n\to 0$ as
$n\to\infty$. Thus the sequence $(\Phi_n)$ is not normal at $0$, and we may apply Zalcman's Lemma (Lemma 1) to it. Replacing $(\Phi_n)$ by a subsequence
if necessary, we thus find that
$$g_n(z)=\Phi_n(z_n+\rho_nz)=1-\frac{2}{w_n}f(\tau_nz_n+\tau_n\rho_nz)\to
g(z)$$
locally uniformly in $\C$, where $|z_n|\leq 1$, $\rho_n>0$, $\rho_n\to 0$,
and $g$ is a non-constant entire function with bounded spherical
derivative. With $\zeta_n=\tau_nz_n$ and $\mu_n=\tau_n\rho_n$ we have
\begin{equation}\label{3c}
g_n(z)=1-\frac{2}{w_n}f(\zeta_n+\mu_nz),\end{equation}
and
\begin{equation}\label{3d}
g^\prime_n(z)=-\frac{2\mu_n}{w_n}f'(\zeta_n+\mu_nz).\end{equation}
We may assume that $\rho_n\leq 1$ and hence $|\zeta_n|\leq\tau_n$ and
$\mu_n\leq\tau_n$ for all $n$.

If $g_n(z)=1$, then $f(\zeta_n+\mu_nz)=0$, hence
$|f'(\zeta_n+\mu_nz)|\leq 1$ by (\ref{3a}). Since $\mu_n\leq\tau_n$, we deduce
that
\begin{equation}\label{3e}
|g'_n(z)|\leq\frac{2\tau_n}{w_n}\quad\mbox{if}\  g_n(z)=1.
\end{equation}
If $g_n(z)=-1$, then $f(\zeta_n+\mu_nz)=w_n$, and hence
$|f^\prime(\zeta_n+\mu_nz)|\leq|w_n|^{1-\varepsilon}$ by our assumption. Thus
\begin{equation}\label{3f}
|g_n^\prime(z)|\leq\frac{2\mu_n}{|w_n|}|w_n|^{1-\varepsilon}\leq\frac{2\tau_n}{
|w_n|^\varepsilon}\quad\mbox{if}\  g_n(z)=-1.
\end{equation}
It follows from the definition of $\tau_n$ that
\begin{equation}\label{3g}
\tau_n=o(|w_n|)^\delta)\quad\mbox{as}\  n\to\infty,
\end{equation}
for any given $\delta>0$.

We deduce from (\ref{3e}), (\ref{3f}) and (\ref{3g}) that $g'(z)=0$
whenever $g(z)=1$ or $g(z)=-1$. Since $g$ has bounded spherical derivative, we 
conclude from
Lemmas~3 and~4 that $g(z)=\cos(az+b).$
Without loss of generality, we may assume that $g(z)=\cos z$ so that $g'(z)=
-\sin z$. In particular, there exist sequences  $(a_n)$ and $(b_n)$ both
tending
to $0$, such that $g_n(a_n)=1$ and $g_n^\prime(b_n)=0$. From (\ref{3e}) we
deduce that
\begin{equation}\label{3h}
|g_n^\prime(a_n)|\leq\frac{2\tau_n}{|w_n|}.
\end{equation}
Noting that $g^{\prime\prime}(z)=-\cos z$ we find that
\begin{equation}
\label{3i}
g_n^\prime(a_n)=g_n^\prime(a_n)-g_n^\prime(b_n)=
\int_{b_n}^{a_n}g^{\prime\prime}_n(z)dz\sim b_n-a_n\end{equation}
as $n\to\infty$, and thus
\begin{equation}\label{3j}
|b_n-a_n|\leq\frac{3\tau_n}{|w_n|}
\end{equation}
for large $n$, by (\ref{3h}). This implies that
\begin{equation}\label{3k}
|g_n(b_n)-1|=|g_n(b_n)-g_n(a_n)|=
\left|\int_{a_n}^{b_n}g^\prime_n(z)dz\right|\leq
2|b_n-a_n|\leq\frac{6\tau_n}{|w_n|}
\end{equation}
for large $n$.

We put
$$h_n(z)=g_n(z+b_n)-g_n(b_n)$$
and note that $h_n(0)=0$, $h^\prime_n(0)=g^\prime_n(b_n)=0$ and
$$h_n(z)\to\cos z-1\quad\mbox{as}\  n\to\infty.$$
It follows that
$$\frac{h_n(z)}{z^2}\to\frac{\cos z-1}{z^2}\quad\mbox{as}\  n\to\infty,$$
which implies that there exists $r>0$ such that
\begin{equation}\label{3l}
\frac{1}{4}\leq\frac{|h_n(z)|}{|z^2|}\leq\frac{3}{4}\quad\mbox{for}
\  |z|\leq r.
\end{equation}
and large~$n$.

Now we fix any $\gamma\in(0,1/2)$ and put
$$c_n=b_n+\frac{1}{|w_n|^\gamma}.$$
Then
$$g_n(c_n)-1=h_n(|w_n|^{-\gamma})+g(b_n)-1$$
and thus, using
(\ref{3k}) and (\ref{3l}) we obtain for large $n$:
\begin{equation} \label{3m}
|g_n(c_n)-1|\leq\left| h_n(|w_n|^{-\gamma})\right|+|g(b_n)-1|
\leq\frac{3}{4|w_n|^{2\gamma}}+\frac{6\tau_n}{|w_n|}
\leq\frac{1}{|w_n|^{2\gamma}}.
\end{equation}
Similarly
\begin{equation}
|g_n(c_n)-1|\geq\left| h_n(|w_n|^{-\gamma})\right|-|g(b_n)-1|
\geq\frac{1}{5|w_n|^{2\gamma}}.\label{3n}
\end{equation}
On the other hand, arguing as in (\ref{3i}), we have
$$g_n^\prime(c_n)=g_n^\prime(c_n)-
g_n^\prime(b_n)=\int_{b_n}^{c_n}g_n^{\prime\prime}(z)dz
\sim b_n-c_n=-\frac{1}{|w_n|^\gamma},$$
and thus
\begin{equation}\label{3o}
|g_n^\prime(c_n)|\geq\frac{1}{2|w_n|^\gamma}
\end{equation} for large $n$.
Put $v_n=\zeta_n+\mu_nc_n.$ Then
$$f(v_n)=\frac{w_n}{2}(1-g_n(c_n))\quad\mbox{and}\quad f^\prime(v_n)=
\frac{w_n}{2\mu_n}g^\prime_n(c_n),$$
by (\ref{3c}) and (\ref{3d}). Hence
\begin{equation}\label{3p}
\frac{1}{10}|w_n|^{1-2\gamma}\leq |f(v_n)|\leq\frac{1}{2}|w_n|^{1-2\gamma},
\end{equation}
by (\ref{3m}) and (\ref{3n}) while
$$|f'(v_n)|\geq\frac{|w_n|^\gamma}{2\mu_n}.$$
Since $|f(v_n)|\geq 1$ for large $n$, by (\ref{3p}),
 this contradicts (\ref{3b}) and (\ref{3g}).

{\em
Mathematisches Seminar, Christian-Albrechts-Universit\"at zu Kiel,

Ludewig-Meyn-Str. 4, 24098 Kiel, Germany

bergweiler@math.uni-kiel.de

\medskip

Department of Mathematics, Purdue University, West Lafayette, IN 47907

eremenko@math.purdue.edu}
\end{document}